\documentclass[preprint]{imsart}

\usepackage{amssymb,amsmath,amsthm}
\usepackage[numbers]{natbib}
\usepackage{hyperref}
\usepackage{cleveref}
\usepackage{enumerate}

\usepackage[T1]{fontenc}
\usepackage{lmodern}
\usepackage{graphicx}
\usepackage{algorithm2e}
\usepackage{subfigure}
\usepackage{ifxetex,ifluatex}
\usepackage{fixltx2e} % provides \textsubscript
\IfFileExists{microtype.sty}{\usepackage{microtype}}{}
\IfFileExists{upquote.sty}{\usepackage{upquote}}{}
\ifnum 0\ifxetex 1\fi\ifluatex 1\fi=0 % if pdftex
  \usepackage[utf8]{inputenc}

\usepackage[a4paper, total={5.2in, 8.2in}]{geometry}

\usepackage{thmtools}
\usepackage{mathtools}

\usepackage{algorithmicx}
\hypersetup{breaklinks=true,
    linkcolor=blue,
    citecolor=magenta,
    colorlinks=true,
pdfborder={0 0 0}}

\usepackage{amsthm}
\usepackage{thmtools}
\usepackage{mathtools}

\usepackage{autonum}

\startlocaldefs

\crefname{thm}{theorem}{theorems}
\crefname{lemma}{lemma}{lemmas}
\crefname{prop}{proposition}{propositions}
\crefname{assumption}{assumption}{assumptions}
\crefname{example}{example}{examples}
\crefname{cor}{corollary}{corollaries}

\declaretheorem[name=Theorem]{thm}
\declaretheorem[name=Proposition,sibling=thm]{prop}

\declaretheorem[name=Assumption]{assumption}

\declaretheorem[name=Example,style=definition]{example}

\newcommand{\E}{\mathbb{E}} 
\newcommand{\R}{\mathbb{R}}

\newcommand{\vxi}{{\boldsymbol{\xi}}}

\newcommand{\proba}[1]{\mathbb{P}\left( #1 \right)}

\newcommand{\psinorm}[1]{\Vert #1 \Vert_{\psi_2}}

\newcommand{\euclidnorm}[1]{\Vert  #1 \Vert_2}
\newcommand{\euclidnorms}[1]{\euclidnorm{ #1 }^2}

\newcommand{\vertiii}[1]{{\vert\kern-0.25ex\vert\kern-0.25ex\vert #1 
            \vert\kern-0.25ex\vert\kern-0.25ex\vert}}

\newcommand{\opnorm}[1]{\vertiii{ #1 }_{2}}

\newcommand{\hsnorm}[1]{\left\Vert #1 \right\Vert_{\mathrm{F}}}

\newcommand{\diag}{\mathrm{diag}}

\newcommand{\CmomentW}{256}
\newcommand{\CmomentWW}{8\sqrt{3}}
\newcommand{\CmomentWWsquared}{192}

\endlocaldefs

\begin{document}

\begin{abstract}
    A concentration result for
    quadratic form of independent subgaussian random variables
    is derived.
    If the moments of the random variables
    satisfy a ``Bernstein condition'',
    then the variance term of the Hanson-Wright inequality
    can be improved.
\end{abstract}

\title{Concentration of quadratic forms under a Bernstein moment assumption}
\runtitle{Concentration of quadratic forms under a Bernstein moment assumption}
\runauthor{P. C. Bellec}

\begin{aug}
    \author{\fnms{Pierre C.} \snm{Bellec}
    }
    \address{Rutgers University, ENSAE and UMR CNRS 9194}
\end{aug}

\date{\today}

\maketitle

\section{Concentration of a quadratic form of subgaussian random variables}

Throughout this note, $A\in\R^{n\times n}$ is a real matrix,
and $\vxi=(\xi_1,...,\xi_n)^T$ is a centered random vector with independent components.
We are interested in the concentration behavior of the random variable
\begin{equation}
    \vxi^T A \vxi - \E[ \vxi^T A \vxi]
    \label{eq:chaos},
\end{equation}
Let $\sigma_i^2 = \E[\xi_i^2]$ for all $i=1,...,n$
and define $D_\sigma = \diag(\sigma_1,...,\sigma_n)$.
If the random variables $\xi_1,...,\xi_n$ are Gaussian,
we have the following concentration inequality.

\begin{prop}[Gaussian chaos of order 2]
    \label{prop:gaussian-chaos}
    Let $\xi_1,...,\xi_n$ be independent zero-mean normal random variables with for all $i=1,...,n$, $\E[ \xi_i^2 ] = \sigma_i^2$.
    Let $A$ be any $n\times n$ real matrix.
    Then for any $x>0$,
    \begin{equation}
        \proba{\vxi^T A \vxi - \E[\vxi^T A \vxi ] > 2 \hsnorm{D_\sigma A D_\sigma} \sqrt{x} + 2 \opnorm{D_\sigma A D_\sigma}x }
        \le
        \exp(-x).
        \label{eq:gaussian-chaos}
    \end{equation}
\end{prop}
A proof of this concentration result can be found in 
\cite[Example 2.12]{boucheron2013concentration}.
We will refer to the term $2 \hsnorm{D_\sigma A D_\sigma} \sqrt{x}$
as the variance term, since if $A$ is diagonal-free,
the random variable $\vxi^TA\vxi$ is centered with variance
\begin{equation}
    \hsnorm{D_\sigma A D_\sigma}^2.
\end{equation}

A similar concentration result is available for subgaussian random variables.
It is known as the Hanson-Wright inequality and is given 
in \Cref{prop:hanson} below.
First versions of this inequality
can be found in 
\citet{hanson1971bound}
and \citet{wright1973bound}, although with a weaker
statement than \Cref{prop:hanson} below since
these results involve $\opnorm{ \, (|a_{ij}|) \, }$ instead of $\opnorm{A}$.
Recent proofs of this concentration inequality
with $\opnorm{A}$ instead of $\opnorm{ \, (|a_{ij}|) \, }$
can be found in \citet{rudelson2013hanson}
or \citet[Theorem A.5]{barthe2013transference}. 
%In the Gaussian case, the quantity
%$\vxi^T A \vxi$ is known as a Gaussian
%chaos of order $2$.
%Concentration bounds for the Gaussian chaos of order $2$
%can be found in \cite[Example 2.12]{boucheron2013concentration}
%or 
%\citet{latala2006estimates}.
%\citet{barthe2013transference} use estimates from 
%\cite{latala2006estimates}
%to prove the result in the subgaussian case.

\begin{prop}[Hanson-Wright inequality \cite{rudelson2013hanson}]
    \label{prop:hanson}
    There exist an absolute constant $c>0$ such that the following holds.
    Let $n\ge 1$ and $\xi_1,...,\xi_n$ be independent zero-mean subgaussian random variables with
    $\max_{i=1,...,n} \psinorm{\xi_i} \le K$ for some real number $K>0$.
    Let $A$ be any $n\times n$ real matrix.
    Then for all $t > 0$,
    \begin{equation}
        \proba{ \vxi^T A \vxi - \E [ \vxi^T A \vxi ]  > t} \le \exp\left( - c \min\left( \frac{t^2}{K^4 \hsnorm{A}^2}, \frac{t}{K^2 \opnorm{A}} \right)\right)
        \label{eq:hanson-t}
    \end{equation}
    where $\vxi=(\xi_1,...,\xi_n)^T$.
    Furthermore, for any $x > 0$, with probability greater than $1-\exp(-x)$,
    \begin{equation}
        \vxi^T A \vxi - \E[\vxi^T A \vxi ] \le c K^2 \opnorm{A}x + c K^2 \hsnorm{A} \sqrt{x}.
        \label{eq:hanson-x}
    \end{equation}
\end{prop}

For some random variables $\xi_1,...,\xi_n$,
the ``variance term'' $K^2 \hsnorm{A} \sqrt{x}$ is far from the variance of the random variable $\vxi^TA\vxi$.
The goal of the present paper is to show that under a mild assumption on the moments of $\xi_1,...,\xi_n$,
it is possible to substantially reduce the variance term.
This assumption is the following.

\begin{assumption}[Bernstein condition on $\xi_1^2,...,\xi_n^2$]
    \label{h:moment}
    Let $K>0$ and assume that $\xi_1,...,\xi_n$
    are independent and satisfy
    \begin{equation}
        \forall p \ge 1, \qquad \E |\xi_i|^{2p} \le \tfrac{1}{2} \; p! \; \sigma_i^2 \; K^{2(p-1)} 
        .
        \label{eq:def-bernstein-rv}
    \end{equation}
\end{assumption}

\begin{example}
    \label{ex:bounded}
    Centered variables almost surely bounded by $K$
    and
    zero-mean Gaussian random variables with variance smaller than $K^2$ satisfy \eqref{eq:def-bernstein-rv}.
\end{example}
\begin{example}[Log-concave random variables] \label{ex:log-concave}
    In \cite{schudy2012concentration},
    the authors consider a slightly stronger condition
    \cite[Definition 1.1]{schudy2012concentration}.
    They consider random variables $Z$ satisfying for any integer $p\ge 1$ and some constant $K$:
    \begin{equation}
        \E[ |Z|^p ] \le \; p \; K \; \E[ |Z|^{p-1} ],
        \label{eq:moment-bounded-schudy}
    \end{equation}
    and they showed in \cite[Section 7]{schudy2012concentration}
    that any distribution that is log-concave
    satisfies (\ref{eq:moment-bounded-schudy}).
    Thus, if $X^2$ is log-concave
    then our assumption (\ref{eq:def-bernstein-rv}) holds.
    See \cite[Section 6]{bagnoli2005log} for a comprehensive
    list of the common log-concave distributions.
\end{example}

The next theorem provides a
concentration inequality for quadratic forms
of independent random variables satisfying
the moment assumption (\ref{eq:def-bernstein-rv}).
It is sharper
than the Hanson-Wright inequality given in \Cref{prop:hanson}.

\begin{thm}
    \label{thm:hanson-moment}
    Assume that the random variable $\vxi=(\xi_1,...,\xi_n)^T$
    satisfies \Cref{h:moment} for some $K>0$.
    Let $A$ be any $n\times n$ real matrix.
    Then for all $t> 0$,
    \begin{equation}
        \proba{ \vxi^T A \vxi - \E[\vxi^T A \vxi]  > t}
        \le
        \exp\left(
            -
            \min\left(
                \frac{t^2}{\CmomentWWsquared K^2 \hsnorm{A D_\sigma}^2},
                \frac{t}{  \CmomentW K^2 \opnorm{A}}
            \right)
        \right),
        \label{eq:hanson-moment-t}
    \end{equation}
    where $D_\sigma = \diag(\sigma_1,...,\sigma_n)$.
    Furthermore, for any $x > 0$, with probability greater than $1- \exp(-x)$,
    \begin{equation}
        \vxi^T A \vxi - \E[\vxi^T A \vxi ]
        \le
        \CmomentW K^2 \opnorm{A} x
        + \CmomentWW K \hsnorm{A D_\sigma} \sqrt{x} 
        .
        \label{eq:hanson-moment-x}
    \end{equation}
\end{thm}

The proof of this result relies on the decoupling
inequality for quadratic forms
\cite{vershynin2011simple} \cite[Theorem 8.11]{foucart2013mathematical}.
% \cite[Theorem 3.1.1]{de1999decoupling}.

If $t$ is small, the right hand side
of (\ref{eq:hanson-moment-t}) becomes
\begin{equation}
    \exp\left( 
        - \frac{t^2}{\CmomentWWsquared K^2 \hsnorm{ A D_\sigma}^2}
    \right),
\end{equation}
whereas the right hand side of the Hanson-Wright inequality (\ref{eq:hanson-t}) becomes
\begin{equation}
    \exp\left( 
        - c \frac{t^2}{K^4 \hsnorm{A}^2}
    \right),
\end{equation}
for some absolute constant $c>0$.
The element of the diagonal matrix $D_\sigma$
are bounded from above by $K$,
so \Cref{thm:hanson-moment}
gives a sharper bound than the Hanson-Wright inequality in this regime.

\section{Proof of \Cref{thm:hanson-moment}}

The goal of this section is to prove \Cref{thm:hanson-moment}.
We start with preliminary calculations that will be useful in the proof.
Let $A$ be any $n\times n$ real matrix.
Let $\lambda > 0$ satisfy
\begin{equation}
    128 \opnorm{A} K^2 \lambda  \le 1,
    \label{eq:condition-lambda}
\end{equation}
and define
\begin{equation}
    \eta = 32  K^2 \lambda^2.
    \label{eq:def-eta}
\end{equation}
The inequality \eqref{eq:condition-lambda} can be rewritten in terms of
$\eta$:
\begin{equation}
    512 K^2 \opnorm{A}^2 \eta \le 1.
    \label{eq:condition-eta}
\end{equation}
Let $A_0$ be the matrix $A$ with the diagonal entries set to $0$.
Then, using the triangle inequality
with $A_0 = A - \diag(a_{11},...,a_{nn})$ 
and $|a_{ii}| \le \opnorm{A}$ for all $i=1,...,n$, we obtain
\begin{equation}
    \opnorm{A_0} \le 2 \opnorm{A}.
    \label{eq:opnorm-Azero}
\end{equation}
Let $B=A_0^T A_0 = (b_{ij})_{i,j=1,...,n}$ and
let $B_0$ be the matrix $B$ with the diagonal entries set to $0$.
Then
\begin{equation}
    \forall i=1,...,n,\qquad
    0 \le b_{ii} = \sum_{j\ne i} a_{ji}^2 \le \opnorm{A}^2.
    \label{eq:bii-positive}
\end{equation}
By using the decomposition $B_0 = B - \diag(b_{11},...,b_{nn})$
and the inequality $\euclidnorms{v+v'}\le 2 \euclidnorms{v} + 2\euclidnorms{v'}$,
\eqref{eq:bii-positive} and \eqref{eq:opnorm-Azero}, we have:
\begin{align}
    \euclidnorms{B_0\vxi} & \le 2 \euclidnorms{B\vxi} + 2 \sum_{i=1}^n b_{ii}^2 \xi_i^2, \\
                          & \le 2 \opnorm{A_0}^2 \euclidnorms{A_0\vxi} + 2 \opnorm{A}^2 \sum_{i=1}^n b_{ii} \xi_i^2, \\
                          & \le 8 \opnorm{A}^2 \euclidnorms{A_0\vxi} + 2 \opnorm{A}^2 \sum_{i=1}^n b_{ii} \xi_i^2.
\end{align}
Combining the previous display with \eqref{eq:condition-eta}, we obtain for any $K>0$:
\begin{align}
    16 K^2 \eta^2 \euclidnorms{B_0\vxi} & 
    \le (512 K^2 \opnorm{A}^2 \eta)\left( \frac{\eta}{4} \euclidnorms{A_0\vxi} + \frac{\eta}{16} \sum_{i=1}^n b_{ii} \xi_i^2\right), \\
    & \le \frac{\eta}{4} \euclidnorms{A_0\vxi} + \frac{\eta}{16} \sum_{i=1}^n b_{ii} \xi_i^2
    .
    \label{eq:preliminary-Azero}
\end{align}

\begin{proof}[Proof of \Cref{thm:hanson-moment}]
    Throughout the proof, 
    let $\lambda>0$ satisfy \eqref{eq:condition-lambda}.
    The value of $\lambda$ will be specified later.

    First we treat the diagonal terms by bounding the moment generating function
    of
    \begin{equation}
        S_{\mathrm{diag}} \coloneqq \sum_{i=1}^n  a_{ii} \xi_i^2 - \sum_{i=1}^n  a_{ii} \sigma_i^2.
    \end{equation}
    Using the independence of $\xi_1,...,\xi_n$ and (\ref{eq:bernstein-condition-1}) with $s=a_{ii} \lambda$ with each $i=1,...,n$:
    \begin{equation}
        \E \exp(\lambda S_{\mathrm{diag}}) \le  \exp\left( \lambda^2 \sum_{i=1}^n a_{ii}^2 \sigma_i^2 K^2 \right),
        \label{eq:bound-mgf-diag}
    \end{equation}
    provided that for all $i=1,...,n$, $2 | a_{ii} | \lambda K^2 \le 1$ which is satisfied as (\ref{eq:condition-lambda}) holds
    and $|a_{ii}| \le \opnorm{A}$.

    Now we bound the moment generating function of the off-diagonal terms.
    Let
    \begin{equation}
        S_{\mathrm{off-diag}} \coloneqq \sum_{i, j=1,...,n : i\ne j} a_{ij} \xi_i \xi_j. 
    \end{equation}
    Let the random vector $\vxi' = (\xi_1',...,\xi_n')^T$ be independent
    of $\vxi$ with the same distribution as $\vxi$.
    We apply the decoupling inequality \cite{vershynin2011simple} (see also \cite[Theorem 8.11]{foucart2013mathematical})
    to the convex function $s\rightarrow \exp(\lambda s)$:
    \begin{equation}
        \E \exp(\lambda S_{\mathrm{off-diag}}) \le \E \exp\left( 4 \lambda \sum_{i,j=1,...,n : i\ne j} a_{ij} \xi_i' \xi_j \right).
    \end{equation}
    Conditionally on $\xi_1,...,\xi_n$,
    for each $i=1,...,n$,
    we use the independence of $\xi_1',...,\xi_n'$
    and (\ref{eq:bernstein-rv-sg}) applied to $\xi_i'$ with $s= 4\sum_{j=1,...,n : i\ne j} a_{ij} \xi_j$:
    \begin{align}
        \E \exp\left( 4 \lambda \sum_{i\ne j} a_{ij} \xi_i' \xi_j \right)
        & \le 
        \E \exp\left( 16  K^2 \lambda^2 \sum_{i=1,...,n} \left( \sum_{j=1,...,n: i\ne j} a_{ij} \xi_j \right)^2 \right), \\
        &=
        \E \exp\left( 16  K^2 \lambda^2 \euclidnorms{A_0 \vxi} \right)
        = 
        \E \exp\left( \frac{\eta}{2} \euclidnorms{A_0 \vxi} \right),
    \end{align}
    where $\eta$ is defined in \eqref{eq:def-eta}
    and $A_0$ is the matrix $A$ with the diagonal entries set to $0$.
    Let $B=A_0^T A_0 = (b_{ij})_{i,j=1,...,n}$.
    Then $\euclidnorms{A_0 \vxi} = \sum_{i=1}^n b_{ii} \xi_i^2 + \sum_{i\ne j} b_{ij} \xi_i \xi_j$.
   
    We use the Cauchy-Schwarz inequality to separate the diagonal terms from the off-diagonal ones:
    \begin{equation}
        \left( \E \exp(\frac{\eta}{2} \euclidnorms{A_0 \vxi}) \right)^2
        \le 
        \E \exp\left( \eta \sum_{i=1}^n b_{ii} \xi_i^2 \right)
        \E \exp\left( \eta \sum_{i\ne j} b_{ij} \xi_i \xi_j \right).
        \label{eq:diag-separation}
    \end{equation}

    For the off-diagonal terms of (\ref{eq:diag-separation}),
    using the decoupling inequality \cite{vershynin2011simple} (see also \cite[Theorem 8.11]{foucart2013mathematical}) we have:
    \begin{equation}
        \E \exp\left(\eta \sum_{i\ne j} b_{ij} \xi_i \xi_j \right)
        \le \E \exp\left( 4 \eta \sum_{i\ne j} b_{ij} \xi_i' \xi_j \right).
    \end{equation}
    Again, conditionally on $\xi_1,...,\xi_n$,
    for each $j=1,...,n$,
    we use (\ref{eq:bernstein-rv-sg}) applied to $\xi_i'$ and the independence of $\xi_1',...,\xi_n'$:
    \begin{align}
        \E \exp\left( 4 \eta \sum_{i\ne j} b_{ij} \xi_i' \xi_j \right)
        & \le 
        \E \exp\left( 16 K^2 \eta^2 \sum_{i=1}^n \left( \sum_{j=1,...,n:\;i\ne j} b_{ij} \xi_j \right)^2 \right), \\
        &=
        \E \exp\left( 16  K^2 \eta^2 \euclidnorms{B_0 \vxi} \right), \\
        & \le
        \E \exp\left( \frac{\eta}{4} \euclidnorms{A_0 \vxi} + \frac{\eta}{16} \sum_{i=1}^n b_{ii} \xi_i^2 \right),
    \end{align}
    where we used the preliminary calculation \eqref{eq:preliminary-Azero} for the last display.
    Finally, the Cauchy-Schwarz inequality yields
    \begin{equation}
        \E \exp\left( 4 \eta \sum_{i\ne j} b_{ij} \xi_i \xi_j' \right)
        \le 
        \sqrt{\E \exp\left( \frac{\eta}{2} \euclidnorms{A_0 \vxi} \right)} 
        \sqrt{\E \exp \left( \frac{\eta}{8} \sum_{i=1}^n b_{ii} \xi_i^2 \right)}.
    \end{equation}
    We plug this upper bound back into \eqref{eq:diag-separation}.
    After rearranging, we find
    \begin{equation}
        \left( \E \exp(\frac{\eta}{2} \euclidnorms{A_0 \vxi}) \right)^{3/2}
        \le 
        \E \exp\left( \eta \sum_{i=1}^n b_{ii} \xi_i^2 \right)
        \sqrt{
            \E \exp \left( \frac{\eta}{8} \sum_{i=1}^n b_{ii} \xi_i^2 \right)
        }.
    \end{equation}
    As $b_{ii}\ge 0$, this implies:
    \begin{equation}
        \E \exp(\frac{\eta}{2} \euclidnorms{A_0 \vxi}) 
        \le 
        \E \exp\left( \eta \sum_{i=1}^n b_{ii} \xi_i^2 \right).
    \end{equation}
    For each $i=1,...,n$, we apply
    (\ref{eq:bernstein-condition-1-bis}) to the variable $\xi_i$
    with $s=b_{ii} \eta \ge 0$.
    Using the independence of $\xi_1^2,...,\xi_n^2$, we obtain:
    \begin{align}
        \E \exp\left( \eta \sum_{i=1}^n b_{ii} \xi_i^2 \right)
        & = \prod_{i=1}^n \E \exp(\eta b_{ii}\xi_i^2), \\
        & \le
        \exp\left(\frac{3}{2} \eta \sum_{i=1}^n b_{ii} \sigma_i^2 \right)
        = 
        \exp\left(\frac{3}{2} \eta \hsnorm{A_0 D_\sigma}^2  \right).
        \label{eq:upper-b-diag}
    \end{align}
    provided that for all $i=1,...,n$, $2K^2 b_{ii} \eta \le 1$
    which is satisfied thanks to \eqref{eq:condition-lambda} 
    and \eqref{eq:bii-positive}.

    We remove $\eta$ from the above displays using its definition
    \eqref{eq:def-eta}:
    \begin{equation}
        \E \exp(\lambda S_{\mathrm{off-diag}})
        \le 
        \exp\left(48 \lambda^2 K^2 \hsnorm{A_0 D_\sigma}^2 \right),
        \label{eq:bound-mgf-off-diag}
    \end{equation}
    where $A_0$ is the matrix $A$ with the diagonal entries set to $0$.

    Now we combine the bound on the moment generating function of
    $S_{\mathrm{diag}}$ and $S_{\mathrm{off-diag}}$,
    given respectively in 
        (\ref{eq:bound-mgf-diag})
        and 
        (\ref{eq:bound-mgf-off-diag}).
    Using the Chernoff bound and the Cauchy-Schwarz inequality:
    we have that for all $\lambda$ satisfying (\ref{eq:condition-lambda}),
    \begin{align}
        \proba{S_{\mathrm{diag}} + S_{\mathrm{off-diag}} >t} 
        & \le
        \exp(- \lambda t) \E [ \exp(\lambda S_{\mathrm{diag}}) \exp(\lambda S_{\mathrm{off-diag}}) ], \\
        & \le
        \exp\left(
            - \lambda t
        \right)
        \sqrt{
            \E [ \exp(2 \lambda S_{\mathrm{diag}}) ]
        }
        \sqrt{
            \E [  \exp(2 \lambda S_{\mathrm{off-diag}}) ]
        }, \\
        & \le
        \exp\left(
            - \lambda t
            + \lambda^2 K^2 \left(
                \sum_{i=1}^n \sigma_i^2 a_{ii}^2 
                + 48 \hsnorm{A_0 D_\sigma}^2 
            \right)
        \right), \\
        & \le 
        \exp\left(
            - \lambda t
            + 48 \lambda^2 K^2 
                \hsnorm{A D_\sigma}^2 
        \right),
        \label{eq:final-rhs}
    \end{align}
    where for the last display we used the equality
    \begin{equation}
        \hsnorm{A D_\sigma}^2  = \sum_{i,j=1,...,n} a_{ij}^2 \sigma_i^2 = \hsnorm{A_0 D_\sigma} + \sum_{i=1}^n a_{ii}^2 \sigma_i^2. 
    \end{equation}
    It now remains to choose the parameter $\lambda$.
    The unconstrained minimum of (\ref{eq:final-rhs})
    is attained at $\bar \lambda = t / (96 K^2 \hsnorm{A D_\sigma}^2)$.
    If $\bar \lambda$ satisfies the constraint
    (\ref{eq:condition-lambda}), then
    \begin{equation}
        \proba{S_{\mathrm{diag}} + S_{\mathrm{off-diag}} >t} 
        \le
        \exp\left(
            \frac{-t^2}{\CmomentWWsquared K^2 \hsnorm{A D_\sigma}^2}
        \right).
    \end{equation}
    On the other hand, if $\bar \lambda$ does not satisfy
    (\ref{eq:condition-lambda}),
    then the constraint 
    (\ref{eq:condition-lambda}) is binding
    and the minimum of 
    (\ref{eq:final-rhs})
    is attained at $\lambda_\mathrm{b} =  1 / ( 128 \opnorm{A} K^2 ) < \bar \lambda$. In this case,
    \begin{equation}
        - t \lambda_\mathrm{b}  
        + \lambda_\mathrm{b}^2 48 K^2 \hsnorm{A D_\sigma}^2
        \le 
        - t \lambda_\mathrm{b}  
        + \lambda_\mathrm{b} \bar \lambda  48 K^2 \hsnorm{A D_\sigma}^2
        =
        - t \lambda_\mathrm{b}  
        + \frac{t}{2} \lambda_\mathrm{b}
        =
        - \frac{t}{\CmomentW K^2 \opnorm{A}}.
    \end{equation}
    Combining the two regimes, we obtain
    \begin{equation}
        \proba{S_{\mathrm{diag}} + S_{\mathrm{off-diag}} >t} 
        \le
        \exp\left(
            -
            \min \left(
                \frac{t^2}{\CmomentWWsquared  K^2 \hsnorm{A D_\sigma}^2},
                \frac{t}{\CmomentW K^2 \opnorm{A}}
            \right)
        \right).
    \end{equation}
    The proof of
    (\ref{eq:hanson-moment-t}) is complete.

    Now we prove
    (\ref{eq:hanson-moment-x}).
    The function
    \begin{equation}
        t \rightarrow  x(t) = 
            \min\left(
                \frac{t^2}{\CmomentWWsquared K^2 \hsnorm{A D_\sigma}^2},
                \frac{t}{\CmomentW K^2 \opnorm{A}}
            \right)
    \end{equation}
    is increasing and bijective from the set of positive real numbers to itself.
    Furthermore, for all $t>0$,
    \begin{equation}
        t \le  \CmomentWW K \hsnorm{A D_\sigma} \sqrt{x(t)} + \CmomentW K^2 \opnorm{A} x(t),
    \end{equation}
    so the variable change $x=x(t)$ completes the proof of 
    (\ref{eq:hanson-moment-x}).
\end{proof}

\section{Technical lemmas: bounds on moment generating functions}

The condition (\ref{eq:def-bernstein-rv}) leads to the following bounds on
the moment generating functions of $X$ and $X^2$,
which are crucial to prove \Cref{thm:hanson-moment}.

\begin{prop}
    \label{prop:def-bernstein-rv}
    Let $K>0$ and let $\xi_i$ be a random variable satisfying 
    (\ref{eq:def-bernstein-rv}) with $\sigma_i^2 = \E[ \xi_i^2 ]$.
    Then for all $s\in\mathbf{R}$:
    \begin{equation}
        \E \exp(s\xi_i) \le \exp(s^2 K^2).
        \label{eq:bernstein-rv-sg}
    \end{equation}
    Furthermore, if $0 \le 2 s K^2 \le 1$, then
    \begin{align}
        \E \exp(s\xi_i^2- s \sigma_i^2) & \le \exp( s^2 \sigma_i^2 K^2), 
        \label{eq:bernstein-condition-1} \\
        \E \exp(s\xi_i^2) & \le \exp\left( \frac{3}{2} s \sigma_i^2\right).
        \label{eq:bernstein-condition-1-bis}
    \end{align}
\end{prop}

Inequality (\ref{eq:bernstein-rv-sg})
shows that a random variable $X$ satisfying 
the moment assumption (\ref{eq:def-bernstein-rv})
is subgaussian and its $\psi_2$ norm is bounded by $K$
up to a multiplicative absolute constant.
The proof of \Cref{prop:def-bernstein-rv} is based on Taylor expansions and some algebra.

\begin{proof}[Proof of \Cref{prop:def-bernstein-rv}]
    To simplify the notation, let $X=\xi_i$ and $\sigma = \sigma_i$.
    We first prove (\ref{eq:bernstein-condition-1}).
    We apply the assumption on the even moments of $X$:
    \begin{align}
        \E \exp(s X^2) & = 1 + s \sigma^2 + \sum_{p\ge 2} \frac{s^p \E X^{2p}}{p!}, \\
                       & \le  1 +  s \sigma^2  + \frac{\sigma^2 s}{2} \sum_{k=1}^\infty (s K^2)^k
        = 1 + s \sigma^2 + \frac{\sigma^2 K^2 s^2 }{ 2 ( 1 - sK^2) },
    \end{align}
    and using the inequality $0< 2s K^2 \le 1$, we obtain:
    \begin{equation}
        \E \exp(s X^2)
        \le 1 + s \sigma^2 + \sigma^2 s^2 K^2 \le \exp(s\sigma^2 + s^2 \sigma^2 K^2),
    \end{equation}
    which completes the proof of 
    (\ref{eq:bernstein-condition-1}).
    Inequality (\ref{eq:bernstein-condition-1-bis}) is a direct consequence of
    (\ref{eq:bernstein-condition-1})
    after applying again the inequality $2 s K^2 \le 1$.

    We now prove (\ref{eq:bernstein-rv-sg}).
    Using the Cauchy-Schwarz inequality and the assumption on the moments for $p=2$,
    we get $\sigma^4 \le \E [ \xi^4 ] \le \sigma^2 K^2$, so $\sigma \le K$.
    Let $p\ge 1$. For the even terms of the expansion of $\E\exp(sX)$, we get:
    \begin{equation}
        \frac{s^{2p} \E X^{2p}}{(2p)!} \le \tfrac{1}{2} (sK)^{2p} \frac{p!}{(2p)!} \le \tfrac{1}{2} \frac{(sK)^{2p}}{p!},
    \end{equation}
    where for the last inequality we used $(p!)^2 \le (2p)!$.
    For the odd terms, by using the Jensen inequality for $p\ge 1$:
    \begin{align}
        \frac{s^{2p+1} \E X^{2p+1}}{(2p+1)!} 
        \le \frac{s^{2p+1} (\E X^{2p+2})^{\frac{2p+1}{2p+2}}}{(2p+1)!} 
        & \le |sK|^{2p+1} \frac{\left(\frac{(p+1)!}{2}\right)^{\frac{2p+1}{2p+2}}}{(2p+1)!}, \\
        & \le \tfrac{1}{2} |sK|^{2p+1} \frac{(p+1)!}{(2p+1)!}.
    \end{align}
    If $|sK| > 1$, we use the inequality $(p+1)!^2 \le (2p+1)!$ to obtain
    \begin{equation}
        \frac{s^{2p+1} \E X^{2p+1}}{(2p+1)!} \le \frac{|sK|^{2(p+1)}}{2((p+1)!)},
    \end{equation}
    and by combining the inequality for the even and the odd terms:
    \begin{align}
        \E \exp(sX) & = 1 + \sum_{p\ge 1} \frac{s^{2p} \E X^{2p}}{(2p)!} + \frac{s^{2p+1} \E X^{2p+1}}{(2p+1)!}, \\
        & \le 1 + \tfrac{1}{2} \sum_{p\ge 1} \frac{(sK)^{2p}}{p!} + \frac{|sK|^{2(p+1)}}{(p+1)!}, \\
        & \le 1 + \sum_{p\ge 1} \frac{(sK)^{2p}}{p!} = \exp(s^2 K^2).
    \end{align}
    If $|sK| \le 1$, we use the inequality $(p+1)! p! \le (2p+1)!$ to obtain
    \begin{equation}
        \frac{s^{2p+1} \E X^{2p+1}}{(2p+1)!} \le \frac{(sK)^{2p}}{2(p!)},
    \end{equation}
    and by combining the inequality for the even and the odd terms:
    \begin{align}
        \E \exp(sX) & = 1 + \sum_{p\ge 1} \frac{s^{2p} \E X^{2p}}{(2p)!} + \frac{s^{2p+1} \E X^{2p+1}}{(2p+1)!}, \\
        & \le 1 + \tfrac{1}{2} \sum_{p\ge 1} \frac{(sK)^{2p}}{p!} + \frac{(sK)^{2p}}{p!}
        = 1 + \sum_{p\ge 1} \frac{(sK)^{2p}}{p!} = \exp(s^2 K^2).
    \end{align}
\end{proof}

\bibliographystyle{plainnat}
\bibliography{../../bibliography/db}

\end{document}